\newcommand{\rem}[1]{}
\documentclass{amsart}
\usepackage{amsfonts,amssymb,amsmath,amsthm,mathrsfs}
\usepackage{url}
\usepackage{graphicx}
\usepackage{epstopdf}
\theoremstyle{definition}

\begin{document}
\author[L.~G.~Molinari]
{Luca~Guido~Molinari}
\address{L.~G.~Molinari: Physics Department Aldo Pontremoli,
Universit\`a degli Studi di Milano and I.N.F.N. sez. Milano,
Via Celoria 16, 20133 Milano, Italy}
\email{Luca.Molinari@unimi.it}

\begin{abstract}  
Several sums of Neumann series with Bessel and trigonometric functions are evaluated, as 
finite sums of trigonometric functions. They arise from a generalization of the Neumann 
expansion of the eigenstates of the Laplacian in regular polygons.
\end{abstract}
\title[Neumann Bessel series]{Some Neumann-Bessel series \\
and the Laplacian on polygons}
\date{9 nov 2020}
\maketitle
\section*{Introduction}
The ground state of the Laplace equation in a regular polygon with Dirichlet boundary conditions at the $n$ sides,
has a natural expression as a Neumann series of Bessel and trigonometric functions,
$$\psi_n (r,\theta)  = J_0(\lambda_n r) +2 {\sum}_{k=1}^\infty h_{k,n} J_{kn} (\lambda_n r) \cos(kn\theta), $$
with coefficients $h_{k,n}$ to be found and eigenvalue $-\lambda_n^2$ that scales with the area. 
For the equilateral triangle and the square, the solutions are known as sums of few trigonometric 
functions of the coordinates $x=r\cos\theta $ and $y=r\sin\theta $. Such solutions 
have a corresponding Neumann expression \cite{Molinari97}. For the square of area $\pi$:
\begin{align}
J_0(\sqrt{2\pi}r)+2{\sum}_{k=1}^\infty J_{4k}(\sqrt{2\pi}r) \cos(4k\theta) = \tfrac{1}{2}\cos (x \sqrt{2\pi}) + \tfrac{1}{2}\cos(y\sqrt{2\pi})
\end{align}
The triangle requires some work to establish the equivalence:
\begin{align}
 J_0(\lambda_3 r) + 2{\sum}_{k=1}^\infty  \frac{\cos(k\pi/2 -\pi/6)}{\cos(\pi/6)} J_{3k}(\lambda_3 r) \cos(3k\theta ) 
 = \tfrac{2}{3\sqrt 3} \sin (\tfrac{4\pi}{3R_3} x +\tfrac{2\pi}{3} ) \\
 -\tfrac{2}{3\sqrt 3} [ \sin [\tfrac{2\pi}{3R_3}(x+y\sqrt 3) -\tfrac{2\pi}{3} ] -\tfrac{2}{3\sqrt 3}
 \sin [\tfrac{2\pi}{3R_3}(x-y\sqrt 3)-\tfrac{2\pi}{3}]   \nonumber 
\end{align}
where, for area $\pi$,  $\lambda^2_3 = \frac{4\pi}{\sqrt 3}$ and $R_3=\tfrac{2}{3}\sqrt{\pi \sqrt 3}$. In \cite{Molinari97}
I also obtained a sum that generalizes the integrable cases $n=3,4$: 
\begin{align}\label{sumn}
f_n(x,y) =& J_0(r) +2 {\sum}_{k=1}^\infty \frac{\cos[nk \tfrac{3\pi}{2}-\tfrac{\pi}{2n}]}{\cos(\tfrac{\pi}{2n})}  
J_{nk}( r) \cos(nk \theta) \\
 =&\frac{1}{n} {\sum}_{\ell=0}^{n-1} \frac{\cos[ r\cos (\theta+\tfrac{2\pi}{n}\ell)+\tfrac{\pi}{2n}] }{\cos (\tfrac{\pi}{2n})}
   \nonumber \\
 =&\frac{1}{n} {\sum}_{\ell=0}^{n-1} \frac{\cos[x \cos (\tfrac{2\pi}{n}\ell) - y \sin (\tfrac{2\pi}{n}\ell)+\tfrac{\pi}{2n}] }{\cos (\tfrac{\pi}{2n})} \nonumber
\end{align} 
For $n\to\infty $ the Riemann sum in the second line is $\int_0^{2\pi} \frac{dt}{2\pi} \cos (r\cos t ) = J_0(r)$; for $n=2$ it is 
$f_2(x,y)=\cos x$.
For $n=6$:
\begin{align}
f_6(x,y)=J_0(r) +2 \sum_{k=1}^\infty (-1)^k
J_{6k}(r) \cos(6k \theta) 
 =\tfrac{1}{3}\cos x +\tfrac{2}{3}\cos(\tfrac{1}{2} x)\cos(\tfrac{\sqrt 3}{2} y) 
\end{align} 
The functions $f_n$ are eigenfunctions of the Laplace operator with eigenvalue $-1$ but, for $n>4$, 
they no longer vanish on the boundary of a $n$-polygon.
\begin{figure}
\begin{center}
\includegraphics*[width=5cm]{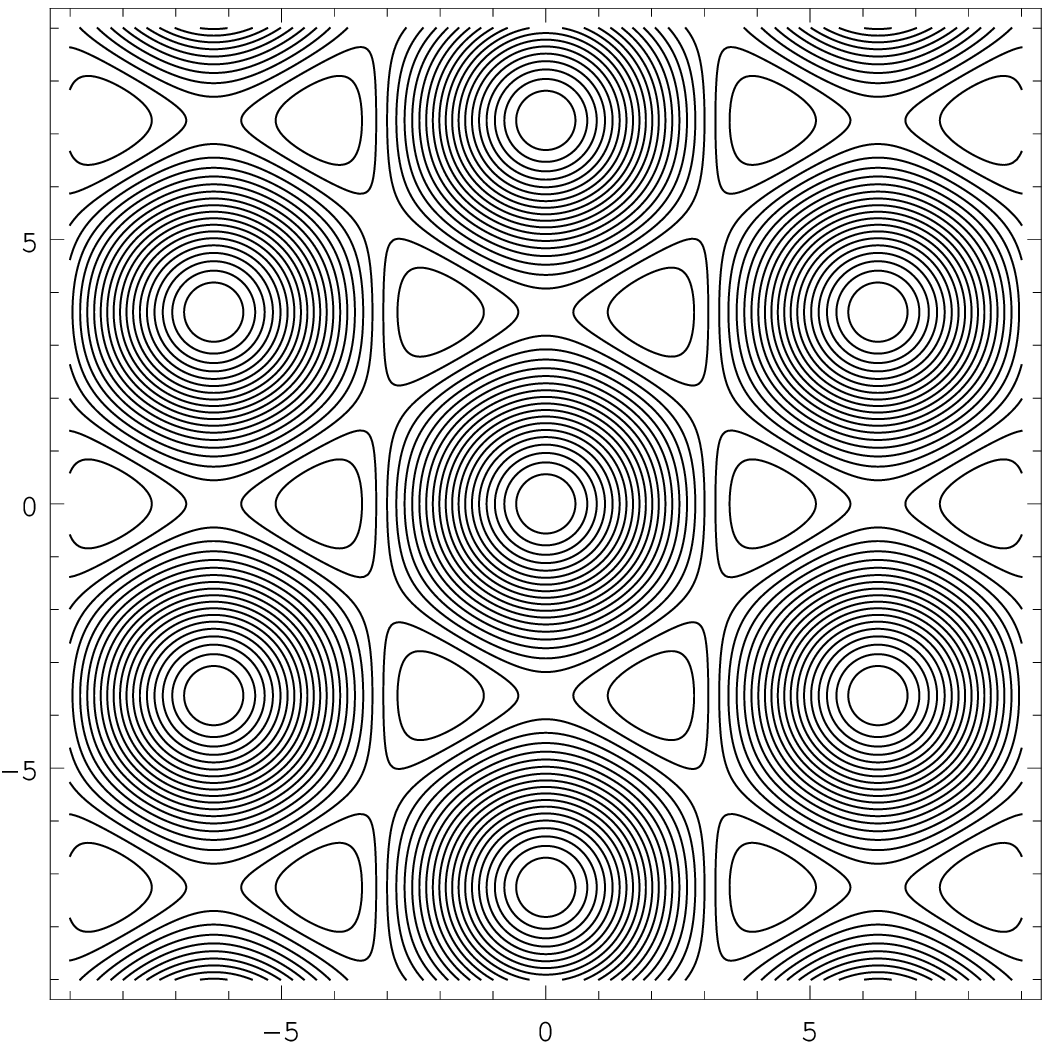}
\includegraphics*[width=5cm]{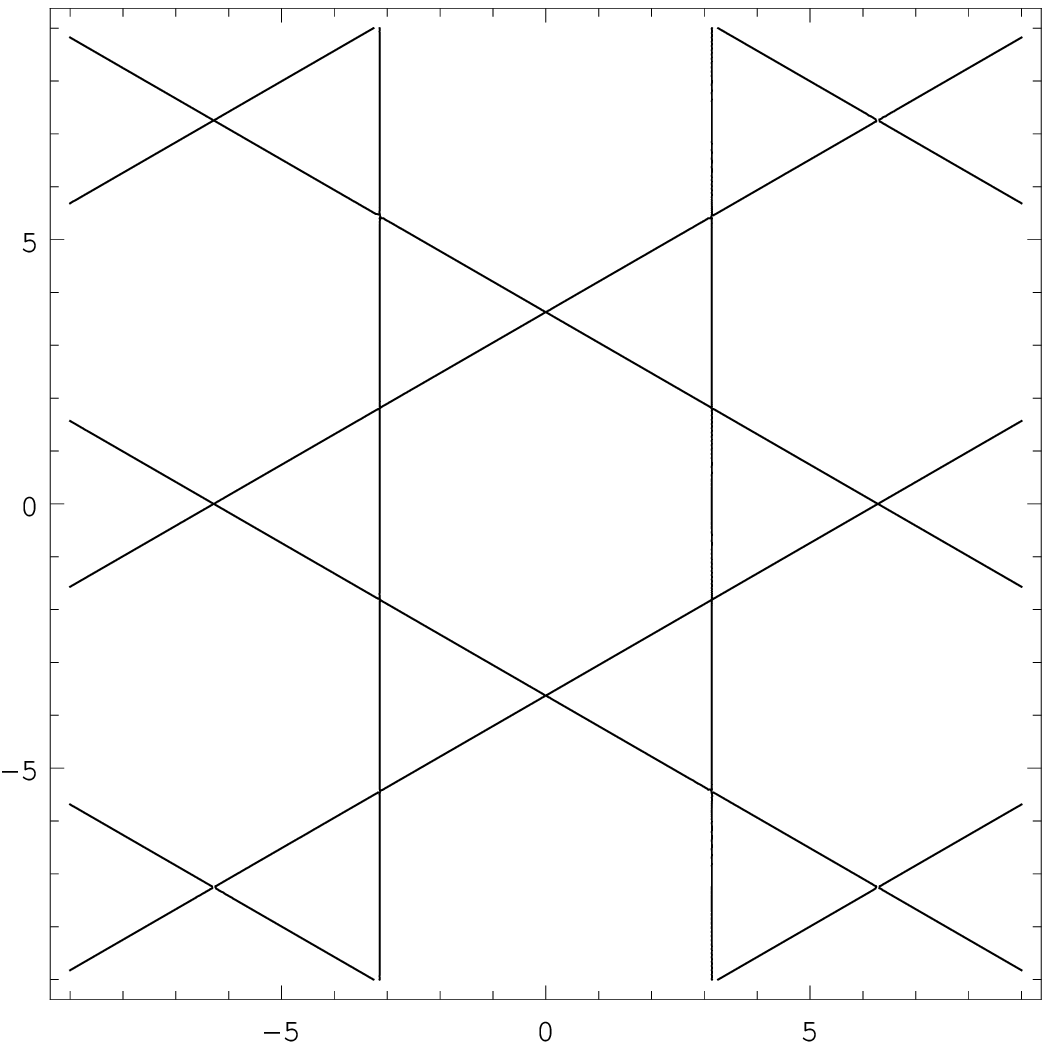}\\
\includegraphics*[width=5cm]{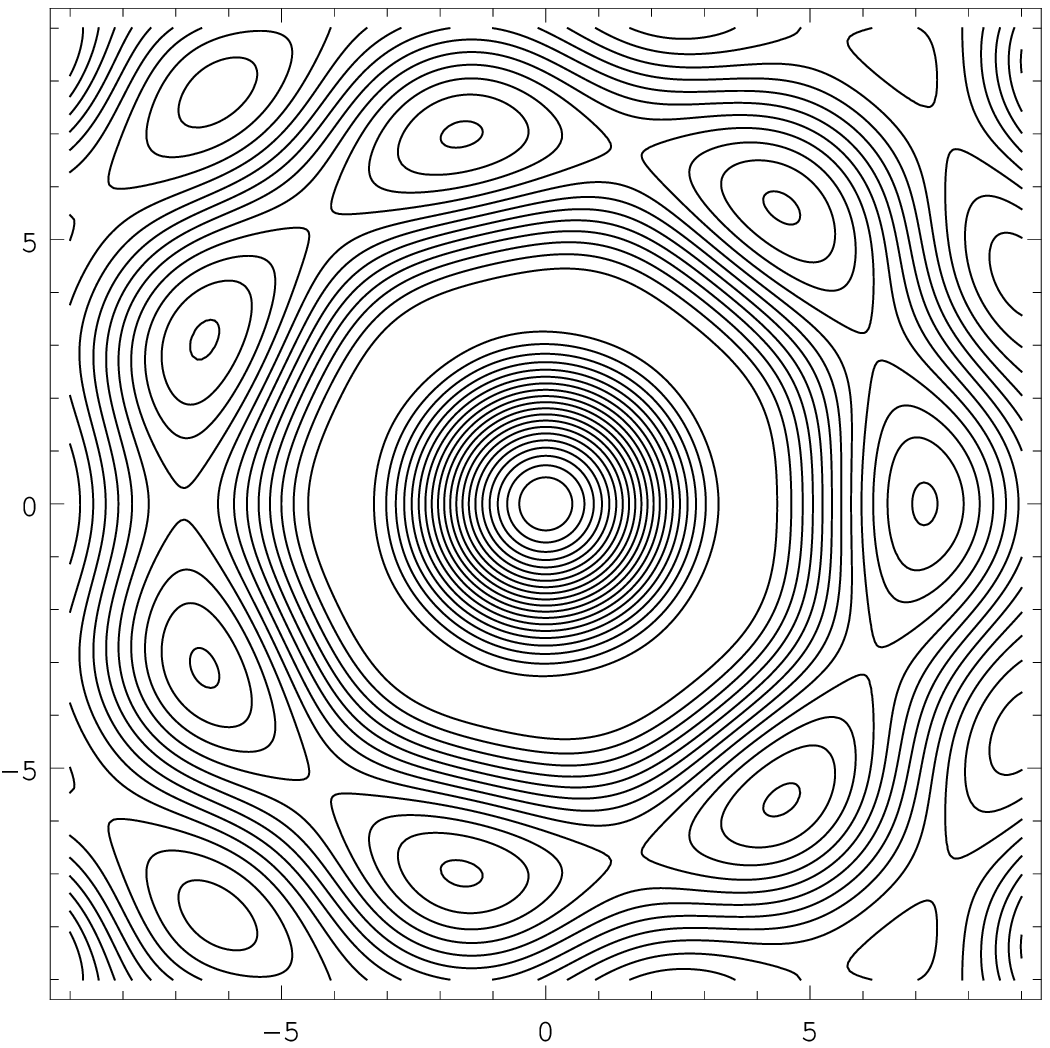}
\includegraphics*[width=5cm]{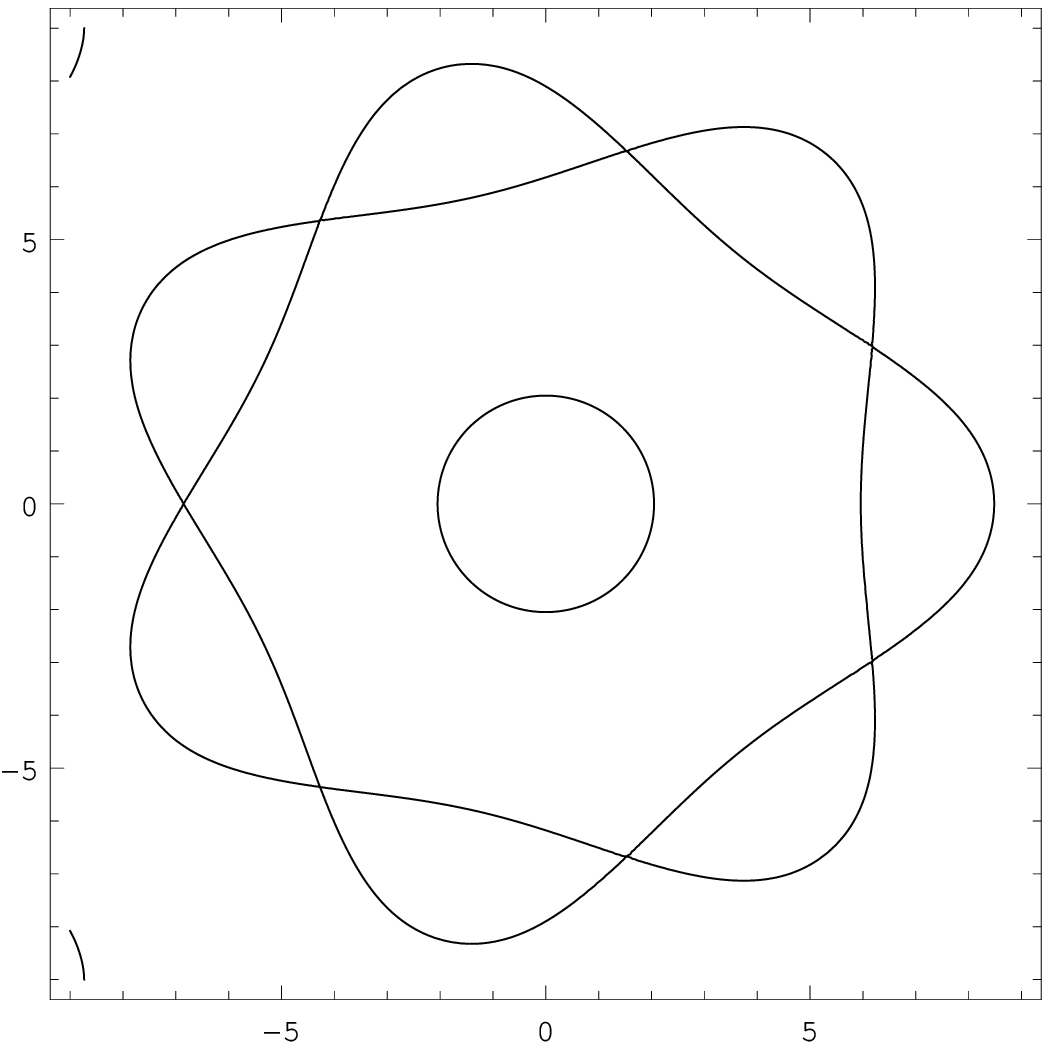}\\
\caption{\label{fig1} 
Left: contour plots of $f_6$ and $f_7$ in eq.\eqref{sumn}. Right: the separatrix line $f_6(x,y)=-1/3$ (Kagom\'e lattice, the eq. 
can be written as $0=\cos(x/2)[\cos(x/2)+\cos(\sqrt 3 y/2)]$. Sites are the doule zeros.),
and  $f_7(x,y)=-1.9633 ...$ (Mathematica)}
\end{center}
\end{figure} 
I only remark that the level curves 
$f_n (x,y)=C_n$ are closed around the origin (where $f_n=1$) up to a separatrix with $n$ self-intersections, with values $C_5=-0.334909$, $C_6=-1/3$, $C_7=0.19633$, etc. The lines are shown in Fig.1.

The Laplacian in polygons has a long history. The ground state beyond the square cannot be finite sums of trigonometric functions, and has been investigated analytically and numerically in $1/n$ expansion (see for example \cite{Molinari97, GS2011,Jones2017}).\\

In this paper I generalize the identity \eqref{sumn}, and obtain a number of new formulas for Neumann series
whose sums contain a finite number of terms. 
For certain values of the parameters, they are identities that are found in the tables 
by Gradshteyn and Ryzhik \cite{GR}, Prudnikov, Brychkov and Marichev \cite{prudnikov}, a recent paper by  Al-Jarrah, Dempsey and Glasser \cite{Glasser}, and two old papers by Takizawa and Kobayasi \cite{Takizawa, Kobayasi}.
In the last ones, the Neumann series appear as correlation functions for the heat flow in coupled harmonic oscillators.

\section*{The summation formula}
The source equation of various sums in this paper is:
\begin{align}
\boxed{ \sum_{k=-\infty}^{+\infty} J_{kn+p}(z) e^{ikny}   
=\frac{1}{n} \sum_{\ell=0}^{n-1} e^{i z\sin (y+\tfrac{2\pi}{n}\ell )-ip (y+\tfrac{2\pi}{n}\ell )} } \label{pippo}
\end{align}
For $y=0$ and even $n$ it is eq.1 in \cite{Takizawa}. Sums of this sort are tabulated only for $n=1,2$ in \cite{prudnikov}.
\begin{proof}
The result follows from the Fourier integral of a Bessel function of integer order.
For fixed $z\in \mathbb C$,  $ \sum_{k=-\infty}^\infty e^{ik\theta} J_{kn+p} (z)$ is uniformly convergent in 
$\theta $ by the bound $|J_{\pm m}(z)|\le C |z/2|^m/m!$ (Nielsen, see \S3.13 in \cite{Watson}).
 \begin{align*}
\sum_{k=-\infty}^\infty e^{ikny} J_{kn+p} (z) =\sum_{k=-\infty}^\infty  e^{ikny} \int_0^{2\pi} \frac{d\theta}{2\pi}  
e^{iz\sin\theta -i(kn+p)\theta}\\
=\sum_{k=-\infty}^\infty  e^{ikny} \sum_{j=0}^{n-1} \int_{\frac{2\pi}{n}j}^{\frac{2\pi}{n}(j+1)} \frac{d\theta}{2\pi}  
e^{iz\sin\theta-ip\theta} e^{-ikn\theta}
\end{align*}
The sums are exchanged:
$= \sum_{j=0}^{n-1}  \sum_{k=-\infty}^\infty  e^{ikny} \int_0^{\frac{2\pi}{n}} \frac{d\theta}{2\pi}  
e^{iz\sin (\theta+\frac{2\pi}{n}j) -ip(\theta+\frac{2\pi}{n}j)}  e^{-ikn\theta} $.
The functions $\sqrt{n/2\pi} \,e^{ikny} $ are 
a complete orthonormal basis in L$^2(0,2\pi/n)$. The infinite sum is the Fourier representation of $\exp[iz\sin (y+\frac{2\pi}{n}j) -ip(y+\frac{2\pi}{n}j)]$.
\end{proof}

%

\subsection*{1}
The case $p=0$ and $y=\frac{\pi}{2}+\alpha $ is an extension with angle $\alpha $ of the 
equations 19 and 20 in \cite{Glasser}, where $\alpha =0$. With $J_{-m}(z)=(-)^mJ_m(z)$:
\begin{align}
J_0(z)+2\sum_{k=1}^{+\infty} e^{i kn\frac{\pi}{2}} J_{kn}(z)  \cos(kn\alpha)
=\frac{1}{n} \sum_{\ell=0}^{n-1} e^{ i z\cos (\alpha +\tfrac{2\pi}{n}\ell )} \label{1p0}
\end{align}
For $n=1$, separation of even and odd parity parts in $z$ gives the Jacobi expansions 
(eqs. 5.7.10.4 and 5 in \cite{prudnikov}):
\begin{align}
J_0(z)+2\sum_{k=1}^{\infty}  (-)^k J_{2k}(z)  \cos(2k\alpha) = \cos (z\cos\alpha)\\
\sum_{k=0}^{\infty} (-)^k J_{2k+1}(z)  \cos[(2k+1)\alpha ]= \tfrac{1}{2}\sin ( z \cos \alpha )
\end{align}
If $\alpha $ is replaced by $\alpha +\pi/2$ they are eqs. 8.514.5 and 6 in \cite{GR} and 10.4, 10.5 in \cite{Kor}:
\begin{align}
J_0(z)+2\sum_{k=1}^{\infty}  J_{2k}(z)  \cos(2k\alpha) = \cos (z\sin\alpha)\\
\sum_{k=0}^{\infty}  J_{2k+1}(z)  \sin[(2k+1)\alpha ]= \tfrac{1}{2}\sin ( z \sin \alpha )
\end{align}
\subsection*{1.1} For $n$ replaced by $2n$, eq.\eqref{1p0} is:
\begin{align}
J_0(z)+2\sum_{k=1}^{\infty}  (-)^{kn} J_{2kn}(z)  \cos(2kn\alpha)
=\frac{1}{2n} \sum_{\ell=0}^{2n-1} \cos [ z\cos (\alpha +\tfrac{\pi}{n}\ell )] \label{1p2}
\end{align}
Since terms $\ell$ and $n+\ell $ are the same, the sum is replaced by
$2\sum_{\ell=0}^{n-1}$.  The value $y=\frac{\pi}{2}$ yields eq.(23) in \cite{Glasser}. \\
For $n=1$ the derivative of \eqref{1p2} in $\alpha =\frac{\pi}{4}$ is:
\begin{align}
2J_2(z)-6J_6(z)+10J_{10}(z) -14J_{14}(z)+ ... = z\tfrac{\sqrt 2}{4}   \sin (z\tfrac{\sqrt 2}{2}) 
\end{align}
For $n=2$ eq.\eqref{1p2} gives
\begin{align}
J_0(z)+2\sum_{k=1}^\infty J_{4k}(z) \cos(4k\alpha) 
=\frac{1}{2} [\cos (z\sin \alpha) +\cos(z\cos \alpha)]\label{13}
\end{align}
The values $\alpha=0, \frac{\pi}{4}$ give eqs. 5.7.1.19. Case $n=3$, $\alpha =0$ gives eq. 5.7.1.21 in \cite{prudnikov}. \\
The derivative of \eqref{13} with $n=4$ is:
\begin{align}
 \sum_{k=1}^\infty k J_{4k}(z) \sin(4k\alpha) =\frac{z}{16} [\sin (z\sin \alpha)\cos \alpha- \sin(z\cos \alpha)\sin \alpha] \label{15}
\end{align}
The expansion in small $\alpha $ gives:
\begin{align}
\sum_{k=1}^\infty k^2 J_{4k}(z)  =\frac{z}{64} (z - \sin z)
\end{align}
\subsection*{1.2} For $n$ replaced by $2n+1$, eq.\eqref{1p0} is:
\begin{align*}
J_0(z)+2\sum_{k=1}^{\infty} e^{i(2n+1)k\frac{\pi}{2}} J_{(2n+1)k}(z)  \cos[(2n+1)k\alpha ]
=\frac{1}{2n+1} \sum_{\ell=0}^{2n} e^{ i z\cos (\alpha +\tfrac{2\pi}{2n+1}\ell )}
\end{align*}
The even-parity and odd-parity parts in the exchange $z\to -z$ are:
\begin{align*}
&J_0(z)+2\sum_{k=1}^{\infty} (-)^k J_{(4n+2)k}(z)  \cos[(4n+2)k\alpha ]
=\frac{1}{2n+1} \sum_{\ell=0}^{2n} \cos[z\cos (\alpha +\tfrac{2\pi}{2n+1}\ell )]\\
&2\sum_{k=0}^{\infty} (-)^{n+k} J_{(2n+1)(2k+1)}(z)  \cos[(2n+1)(2k+1)\alpha ]
=\frac{1}{2n+1} \sum_{\ell=0}^{2n} \sin [z\cos (\alpha +\tfrac{2\pi}{2n+1}\ell )]
\end{align*}
Examples of the second equation are
\begin{align}
&\sum_{k=0}^{\infty} (-)^k J_{6k+3}(z)  \cos[(6k+3)\alpha ]
=-\frac{1}{6} \sum_{\ell=0}^{2} \sin [z\cos (\alpha +\tfrac{2\pi}{3}\ell )] \label{hexagon_triangle}\\
&\sum_{k=0}^{\infty} (-)^k J_{10k+5}(z)  \cos[(10k+5)\alpha ]
=\frac{1}{10} \sum_{\ell=0}^{4} \sin [z\cos (\alpha +\tfrac{2\pi}{5}\ell )] \label{deca}
\end{align}
The first equation with $\alpha =\pi$ is eq.22 in \cite{Glasser}.\\ 
Both sums are eigenfunctions of the Laplacian
with eigenvalue $\lambda=-1$ (see Fig. 2). The sum \eqref{deca}, with $z=r$, $x=r\cos\alpha $ and $y=r\sin\alpha $, is
\begin{align}
f(x,y) =
 \sin x -2\sin (x \cos\tfrac{\pi}{5})\cos (y  \sin\tfrac{\pi}{5}) +2\sin (x \cos\tfrac{2\pi}{5})\cos (y  \sin\tfrac{2\pi}{5})
\label{DECX}
\end{align}

\begin{figure}
\begin{center}
\includegraphics*[width=4.5cm]{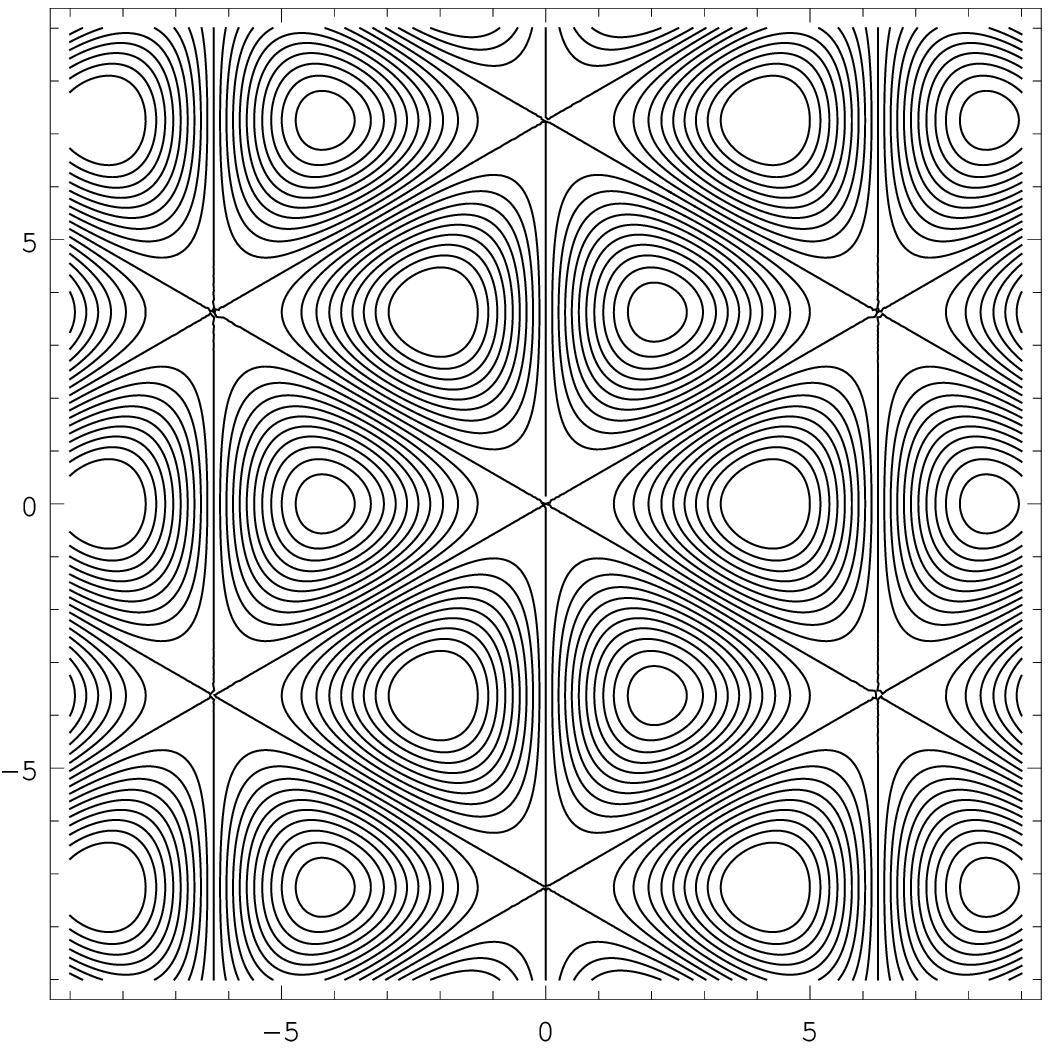}
\includegraphics*[width=4.5cm]{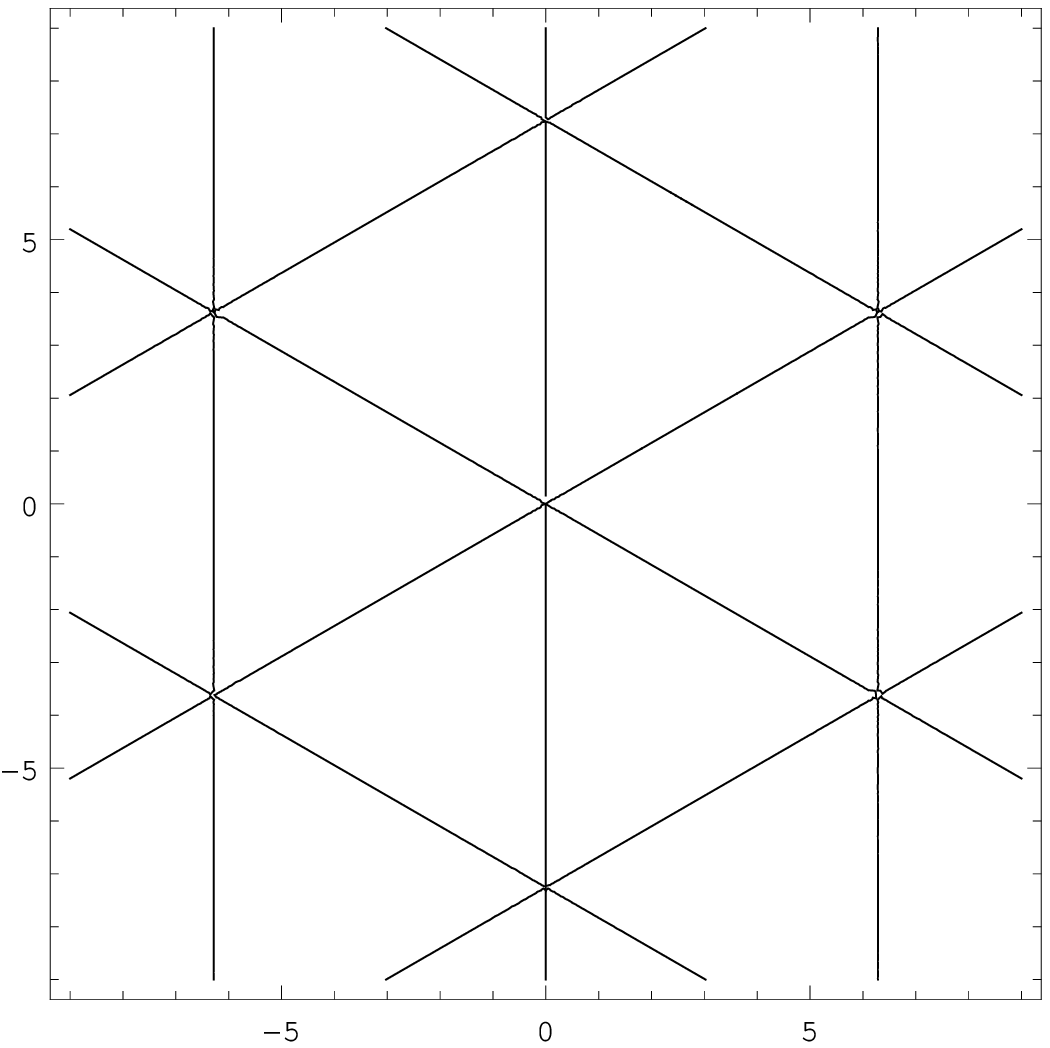}\\
\includegraphics*[width=4.5cm]{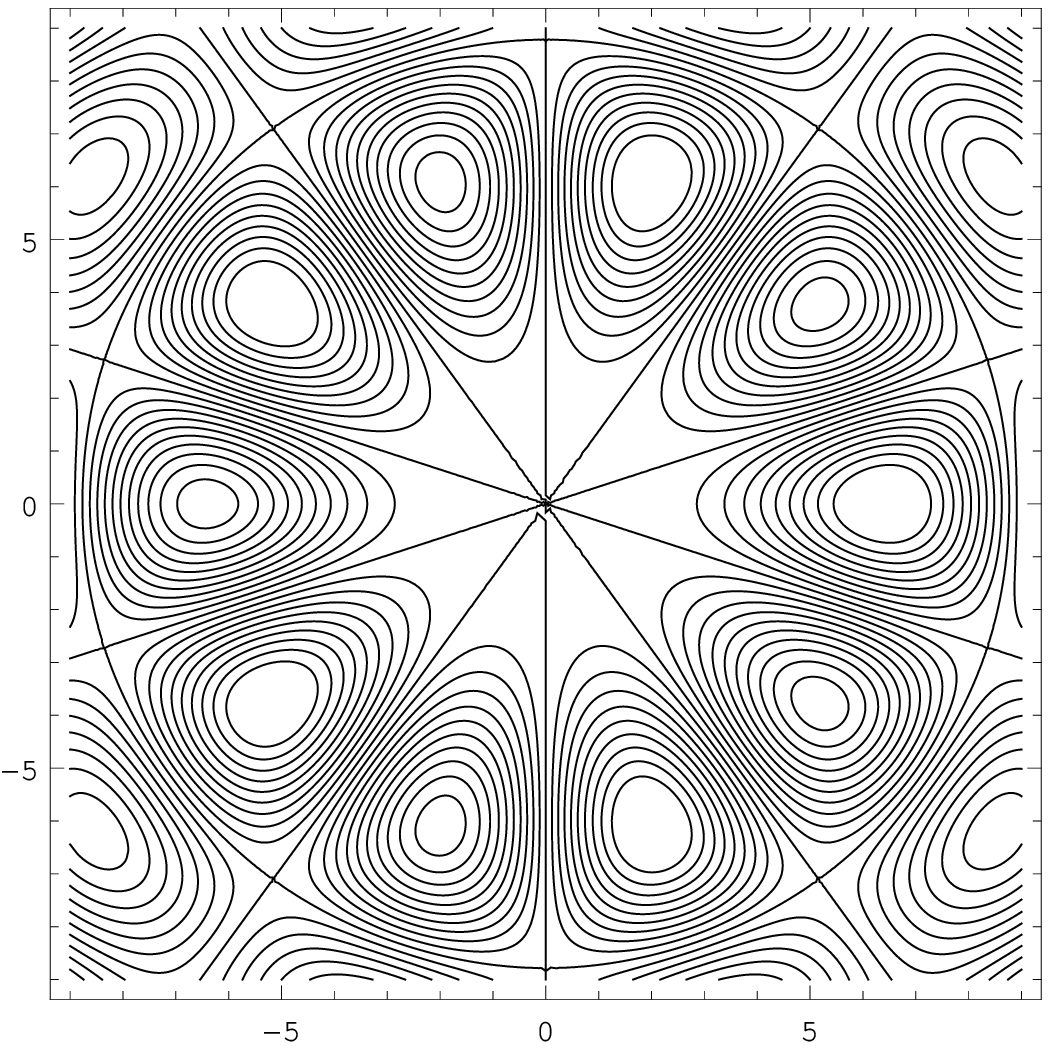}
\includegraphics*[width=4.5cm]{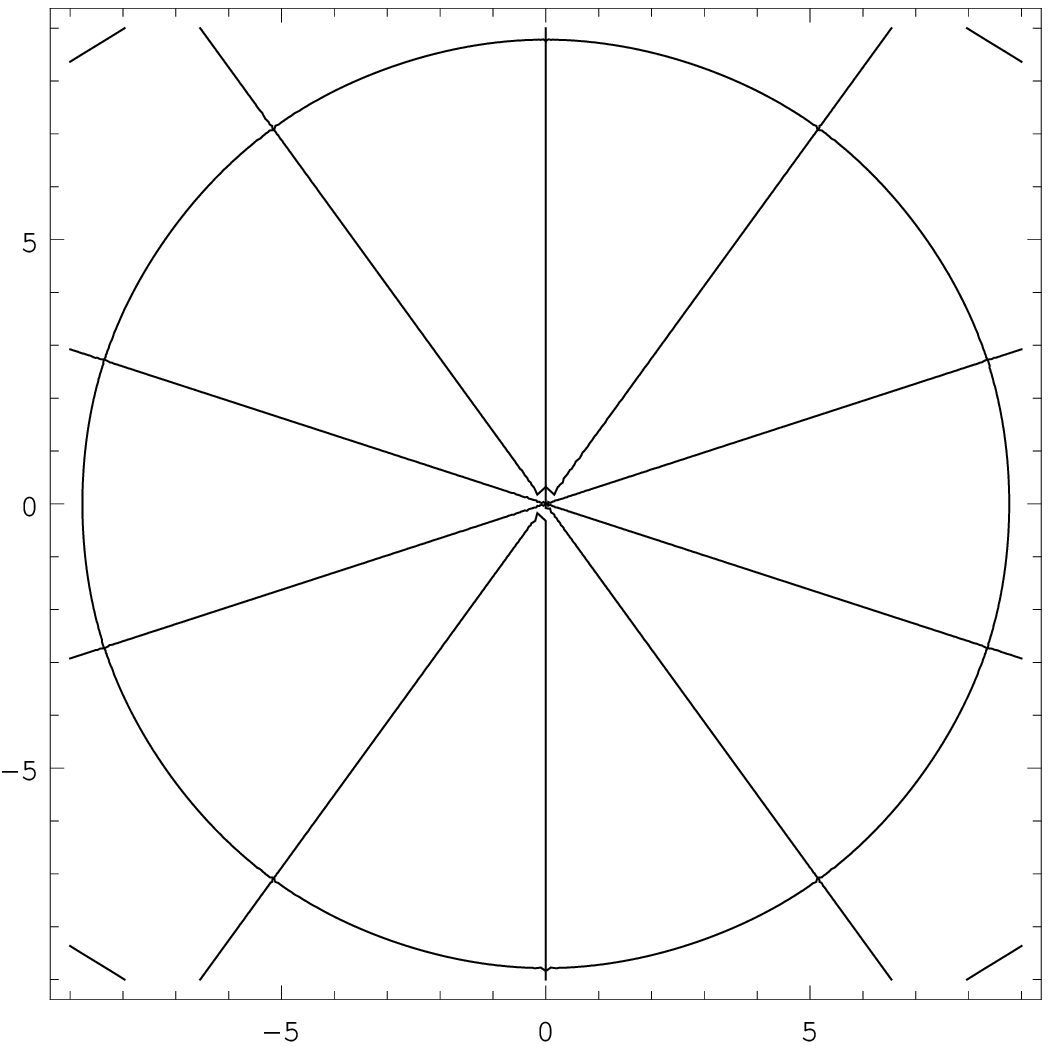}
\caption{Contour plots of the sums \eqref{hexagon_triangle} and \eqref{DECX}. The first is the ground state of the equilateral triangle (no nodal lines) and an excited state of the hexagon. The second function is zero on a line close to the first zero of $J_5(r)$.}
\end{center}
\end{figure} 

\subsection*{1.3}
In \eqref{pippo} with $p=0$, multiply by $\exp(i\beta)$ ($\beta $ real) and take the
real part. The left hand side becomes:
\begin{align*}
&J_0(x)\cos\beta +{\sum}_{k=1}^\infty  J_{kn}(x) \; {\rm Re}[e^{i\beta} (e^{ikny} + e^{-ikn (y+\pi) })]\\
&=J_0(x)\cos\beta +2 {\sum}_{k=1}^\infty  \cos(\beta - kn\tfrac{\pi}{2}) J_{kn}(x) \cos[kn (y+\tfrac{\pi}{2})]
\end{align*} 
The identity \eqref{sumn} is obtained, when $\beta=\tfrac{\pi}{2n}$ and $y +\tfrac{\pi}{2}=\theta+\pi$. 

\subsection*{2} Parseval's identity is applied to \eqref{pippo}:
\begin{align*}
\sum_{k\in\mathbb Z} J^2_{kn+p}(x) &= \frac{1}{n^2}\sum_{k\ell} e^{ip \tfrac{2\pi}{n}(k-\ell) }
\int_0^{2\pi} \frac{dy}{2\pi} e^{ix \sin (y-\tfrac{\pi}{n}(k-\ell)) - ix\sin (y+\tfrac{\pi}{n}(k-\ell)) }\nonumber\\
&= \frac{1}{n^2}\sum_{k\ell} e^{ip \tfrac{2\pi}{n}(k-\ell) }
\int_0^{2\pi} \frac{dy}{2\pi} e^{-i2x \cos y \sin(\tfrac{\pi}{n}(k-\ell))  } \nonumber\\
&= \frac{1}{n^2}\sum_{k\ell} e^{ip \tfrac{2\pi}{n}(k-\ell) }
J_0[2x \sin(\tfrac{\pi}{n}(k-\ell))]\nonumber \\
&=\frac{1}{n}+\frac{2}{n^2}\sum_{k=1}^{n-1} k \cos(\tfrac{2\pi}{n}kp) J_0(2x \sin \tfrac{\pi k}{n} )
\end{align*}
The last sum is unchanged if $k$ is replaced by $n-k$: 
\begin{align}
\sum_{k\in\mathbb Z} J^2_{kn+p}(x) 
=\frac{1}{n}+\frac{1}{n}\sum_{k=1}^{n-1} \cos(\tfrac{2\pi}{n}kp) J_0(2x \sin \tfrac{\pi k}{n} ) \label{Parseval0}
\end{align}
The left-hand side is 
$J_p(x)^2+\sum_{k=1}^\infty J^2_{kn+p}(x)+J^2_{kn-p}(z ) $. \\
For the special case $p=0$ and $n\to 2n$ in \eqref{Parseval0}, the sum is amenable to eq.29 in \cite{Glasser}:
\begin{align}
J_0^2(x)+  2\sum_{k=1}^{\infty}  J_{2kn}^2 (x)= \frac{1}{2n}+ \frac{1}{2n} J_0(2x) + \frac{1}{n} \sum_{k=1}^{n-1}  J_0 (2x\cos \tfrac{\pi}{2n}k ) \label{Parseval1}
\end{align}

\subsection*{2.1}
If $n\to 2n$ and $p=n$ in \eqref{Parseval0}, with simple steps one obtains:
\begin{align}
\sum_{k=0}^\infty J^2_{(2k+1)n}(x) = \frac{1}{4n}+ \frac{(-)^n}{4n} J_0(2x) +\frac{1}{2n}\sum_{\ell=1}^{n-1} (-1)^\ell 
J_0(2x \sin\tfrac{\pi\ell}{2n})  \label{Parseval2}
\end{align}

\subsection*{3} In eq.\eqref{pippo} the variable $y$ is shifted to $y+2t$. The equation is multiplied by $e^{iz'\sin y - iqy}$ and integrated in $y$:
\begin{align}
\sum_{k=-\infty}^{+\infty} J_{p+kn}(z)& J_{q-kn} (z') e^{i(kn+p)2t} \nonumber \\
=&\frac{1}{n} \sum_{\ell=0}^{n-1}e^{-ip\tfrac{2\pi}{n}\ell} 
 \int_0^{2\pi} \frac{dy}{2\pi}e^{i z\sin (y+2t+\tfrac{2\pi}{n}\ell )+iz'\sin y -i(p+q)y} \label{pippo2} 
\end{align}
In the integral, the shift $y$ to $y-t-\frac{\pi}{n}\ell $ changes the exponent to
$$ i (z+z')\sin y \cos(t+\tfrac{\pi}{n}\ell )+i(z-z')\cos y\sin (t+\tfrac{\pi}{n}\ell ) -i(p+q)(y-t-\tfrac{\pi}{n}\ell ) $$
\subsection*{3.1} With $z=z'$ we obtain eq.1 in \cite{Kobayasi}:
\begin{align}
\sum_{k=-\infty}^{+\infty} J_{p+kn}(z) J_{q-kn} (z) e^{2iknt }  
=\frac{1}{n} \sum_{\ell=0}^{n-1}e^{-i(p-q)(t+\tfrac{\pi}{n}\ell )}
  J_{p+q}[2z \cos(t+\tfrac{\pi}{n}\ell ) ] 
\end{align}
For $n=1,2$ it is (with a shift of the index $k$ in the first identity and renaming of parameter):  
\begin{align}
&\sum_{k=-\infty}^{+\infty} J_{k}(z) J_{p-k} (z) e^{2ikt }  
=e^{ipt} J_p(2z \cos t )\\
&\sum_{k=-\infty}^{+\infty} J_{p+2k}(z) J_{q-2k} (z) e^{4ikt }  
=\tfrac{1}{2}e^{-i(p-q)t} [J_{p+q}(2z \cos t )+ i^{p-q} J_{p+q}(2z \sin t)]
\end{align}
The first one is eq.8.530 \cite{GR}. The second one, for $t=\frac{\pi}{4}, \frac{\pi}{2}$, becomes:
\begin{align}
&\sum_{k=-\infty}^{+\infty} (-)^kJ_{p+2k}(z) J_{q-2k} (z) 
= J_{p+q}(z \sqrt 2 ) \cos [(p-q)\tfrac{\pi}{4}]\\
&\sum_{k=-\infty}^{+\infty} J_{p+2k}(z) J_{q-2k} (z) = \tfrac{1}{2} J_{p+q}(2z)
\end{align}
For $p=q$ the first one is eq. 5.7.11.25 \cite{prudnikov}.
\subsection*{3.2} Eq.\eqref{pippo2} with $p=q=0$ and $t=0$ is:
\begin{align}
J_0(z)J_0(z')+2\sum_{k=1}^\infty (-)^{kn} J_{kn}(z) J_{kn}(z')
=\frac{1}{n} \sum_{\ell=0}^{n-1} \int_0^{2\pi}\frac{dy}{2\pi} e^{i z\sin (y+\tfrac{2\pi}{n}\ell )+iz'\sin y} \nonumber
\end{align}
For $n=1, 2$ they are eqs. 5.7.11.1 and 5.7.11.18 in \cite{prudnikov} and, for $z=z'$: eqs. 31, 32 in \cite{Glasser}.
A new example is:
\begin{align}
\sum_{k=1}^\infty  J_{4k}(x) J_{4k}(y) =\tfrac{1}{8} [J_0(x+y) +J_0(x-y)- 4J_0(x)J_0(y)+ 2J_0 (\sqrt{x^2+y^2})]
\end{align}
\rem{+++++++
\subsection*{3.3} Eq.\eqref{pippo2} with $p=q$ and $n=2p$ is:
\begin{align}
J_p(x)J_p(y)+&\sum_{k=1}^\infty (-1)^p[J_{2kp+p}(x) J_{2kp-p}(y)+ J_{2kp-p}(x)J_{2kp+p}(y)] \\
&=\frac{1}{2p} \sum_{\ell=0}^{2p-1}  \int_0^{2\pi} \frac{ds}{2\pi}e^{i (x+y)\sin s\cos(\tfrac{\pi}{2p}\ell )+i(x-y)\cos x\sin(\tfrac{\pi}{2p}\ell) -2ips} \nonumber\\
&=\frac{1}{2p} \sum_{\ell=0}^{2p-1}  \int_0^{2\pi} \frac{ds}{2\pi}e^{i R_\ell \sin (s+\theta_\ell) -2ips} \nonumber\\
&=\frac{1}{2p} \sum_{\ell=0}^{2p-1} \frac{J_{2p}(R_\ell )}{R_\ell^{2p}} [(x+y)\cos(\tfrac{\pi}{2p}\ell) -
i (x-y)\sin(\tfrac{\pi}{2p}\ell) ]^{2p}
\end{align}
where $R_\ell^2 = x^2+y^2 +2xy \cos(\tfrac{\pi}{p}\ell)$, $(x+y)\cos(\tfrac{\pi}{2p}\ell) =R_\ell \cos\theta_\ell$, 
$(x-y)\sin(\frac{\pi}{2p}\ell) =R_\ell \sin\theta_\ell$.
+++++++}
\subsection*{4} Multiplication of \eqref{pippo} by $\exp (-a y)$ ($a>0$) with $p=0$ and $n=1$,  and 
integration on $\mathbb R^+$ give:
\begin{align}
\frac{1}{a}J_0(z)+
\sum_{k=1}^\infty J_{2k}(z) \frac{2a}{a^2+4k^2}  +
\sum_{k=0}^\infty J_{2k+1}(z) \frac{2i(2k+1)}{a^2+(2k+1)^2}  
= \int_0^\infty dy\, e^{i z\sin y-ay} \nonumber
\end{align}
The integral in the right-hand side is done by series expansion, with eqs.3.895.1 and 3.895.4 \cite{GR}. The even and odd
terms are: 
\begin{align}
&\frac{1}{a} J_0 (z) + \sum_{k=1}^\infty J_{2k}(z) \frac{2a}{a^2+4k^2} 
=\sum_{k=0}^\infty (-)^k \frac{z^{2k}}{a(a^2+4)...(a^2+4k^2)}\\
&\sum_{k=0}^\infty J_{2k+1}(z) \frac{2(2k+1)}{a^2+(2k+1)^2}  
= \sum_{k=0}^\infty (-)^k \frac{z^{2k+1}}{(a^2+1)(a^2+9)...(a^2+(2k+1)^2)}
\end{align}
More and more identities can be obtained by derivation, or integration with functions. Here I limited myself to simple and,
hopefully, useful examples.

\end{document}